\documentclass[journal]{IEEEtran}
\usepackage{amssymb}
\usepackage{amsthm}
\usepackage[cmex10]{amsmath}
\usepackage{graphicx}
\usepackage{colortbl}
\usepackage{cite}

\newtheorem{theorem}{Theorem}

\newtheorem{remark}{Remark}
\newtheorem{lemma}{Lemma}

\ifCLASSINFOpdf
\else
\fi
\begin{document}
%
\title{Output Feedback Control for Irregular LQ Problem}
%
%
%

\author{Juanjuan~Xu, Huanshui~Zhang
\thanks{This work was supported by the National Natural Science Foundation of China (Grant
Nos. 61633014, 61573221, and 61873332) and the Qilu Youth Scholar Discipline Construction Funding from
Shandong University.}
\thanks{J. Xu and H. Zhang are with School of Control
Science and Engineering, Shandong University, Jinan, Shandong, P.R.China 250061.
        (e-mail: juanjuanxu@sdu.edu.cn, hszhang@sdu.edu.cn)}
}

\maketitle

\begin{abstract}
In this paper, we study the irregular output feedback linear quadratic (LQ) control problem, which is a continuous work of previous works for irregular LQ control \cite{zhangSCIS} where the state is assumed to be exactly known priori.  Different from the classic output feedback LQ control problem, the cost function must be modified to guarantee the solvability.
In the framework of the modified cost function, it is shown that the separation principle holds and the explicitly optimal controller is given in the feedback form of the Kalman filtering.
In particular, the feedback gain is calculated by two Riccati equations independently of the Kalman filtering. The key technique is the ``two-layer optimization" approach.
We also emphasize that the optimal controller at the terminal time is required to be deterministic.
\end{abstract}


%
\IEEEpeerreviewmaketitle

\section{Introduction}\label{s1}

Linear quadratic (LQ) optimal control problem is one of the most fundamental problems~\cite{hszhang}. It was first studied by R. Bellman in 1958 and
solved in a linear state feedback control form by R. E. Kalman in 1960's. Since then, extensive research has been paid to LQ control problem~\cite{Letov,anderson,lewis}.
Initially, the controller's weight matrix is assumed to be positive-definite. The unique solution is thus explicitly given in terms of a standard Riccati equation~\cite{Kalman}.
However, through continuous research on engineering systems, economic models and natural resource production,
it is found that the assumption of positive-definite weight matrix is too strict.
In fact, the case of singular weight matrix is widespread in which the LQ problem is termed as  singular control~\cite{Bell}.
A typical example occurs in aerospace applications whose aim is to control rockets to minimize the propellant expenditure.

In view of its broad applications, singular control problem has also been widely studied. For instance, \cite{Ho} studied the case of zero weighting matrix of the control in the cost function. It is obtained that
the problem is solvable for some specific initial values and an impulse control must be applied to guarantee the solvability for arbitrary initial values. For the general singular case,
some sufficient conditions have been given by using the transformation approach~\cite{Gurman,Moore,Williems}. In~\cite{Gabasov,Krener, Hoehener, Bonnans, XuZhang}, singular control problems were studied by
exploring the higher order maximum principle. Noting that this approach fails when the higher derivatives vanish~\cite{Gabasov}.
By adding a perturbation in the singular weighting matrix, the problem was solved by applying the standard LQ control~\cite{chenhanfu}. With the perturbation approach,
\cite{sunliyong} derived the open-loop solution of irregular LQ problem by applying a minimize convergence sequence and gave the closed-loop controller under the regularity of the Riccati equation. Though the problem has been widely studied from 1950's, some fundamental problems remain to be solved. Until recently, \cite{zhangSCIS} obtained the complete solution for the irregular LQ control with a new approach of ``two-layer optimization". It has been shown that the essential difficulty to solve the singular
control problem lies in the irregularity.
Moreover, the irregular LQ control is totally different from the regular one due to that the irregular controller must guarantee the terminal state constraint of $P_1(T)x(T)=0.$

In the above mentioned work, there is a presupposition that the state information is perfectly available to designing the controller.
However, in practical applications, 
the state information is usually  affected by stochastic noises when transforming to the controller.
This motivates us to further study the irregular output feedback control problem.

It should be pointed out that the iregular LQ problem becomes much more involved when the system is affected by additive noises.
In order to illustrate this point more clearly, we present the following example for LQG.

Consider the system
\begin{eqnarray}
dx(t)&=&u(t)dt+dw(t),~x(0)=x_0,\label{1.1}
\end{eqnarray}
where $w(t)$ is a scalar standard Brownian motion. The LQG problem is to seek $u(t)$ minimizing the following cost function:
\begin{eqnarray}
J_1=E[x^2(T)]. \label{1.2}
\end{eqnarray}
The related LQ problem of the above LQG is as
\begin{eqnarray}
\min_{u}x^2(T)~~\mbox{s.t.}~~\dot{x}(t)=u(t). \label{1.3}
\end{eqnarray}
By \cite{zhangSCIS}, it is easy to know that LQ problem (\ref{1.3}) is irregular,  and the open-loop optimal control is as $u(t)=-\frac{x_0}{T}$ and
the closed-loop optimal solution is $u(t)=\frac{x(t)}{t-T},$  the optimal cost is $0.$
%

However, it is easy to verify that the LQG problem is unsolvable with the cost function (\ref{1.2}). In fact, with the aid of the stochastic maximum principle, the optimal solution satisfies the following forward and backward stochastic differential equations (FBDEs):
\begin{eqnarray}
dx(t)&=&u(t)dt+dw(t),~x(0)=x_0,\nonumber\\
dp(t)&=&q(t)dw(t),~p(T)=x(T),\nonumber\\
0&=&p(t).\nonumber
\end{eqnarray}
Thus, it must hold almost surely from $p(t)=0$ that
\begin{eqnarray}
x(T)=0.\label{op13}
\end{eqnarray}
However, by solving the stochastic differential equation $dx(t)=u(t)dt+dw(t),$ we have
\begin{eqnarray}
x(T)=x_0+\int_0^Tu(s)ds+\int_0^Tdw(s).\nonumber
\end{eqnarray}
Let $u(t)$ be open-loop, it gives that
\begin{eqnarray}
E\|x(T)\|^2&=&E\|x_0+\int_0^Tu(s)ds\|^2+E\|\int_0^Tdw(s)\|^2\nonumber\\
&=&E\|x_0+\int_0^Tu(s)ds\|^2+T\nonumber\\
&\geq&T.\label{op14}
\end{eqnarray}
This contradicts the solvability condition (\ref{op13}). Thus, the LQG problem is unsolvable.\\

To guarantee the solvability of the stochastic control problem, we modify the cost function as
\begin{eqnarray}
J_2=[Ex(T)]^2, \nonumber
\end{eqnarray}
Then the modified LQG problem is
$$\min_{u} J_2~~\mbox{s.t.}~~dx(t)=u(t)dt+dw(t).$$
The open-loop and closed-loop optimal controllers are $u(t)=-\frac{x_0}{T}$ and
$u(t)=\frac{x(t)}{t-T}$ respectively, and the corresponding optimal cost is $0$.  This is consistent with the deterministic case.

Based on the discussions of the above example, it is seen that the irregular control problem with additive noises becomes much more complex.
In this paper, we will study the output feedback control for the irregular LQ problem.
The optimal controller is given by considering a modified cost function as in the above. In the frame of the new cost function, it is shown that
the separation principle holds. To be specific, the optimal control is in the feedback form of the Kalman filtering.
By applying the ``two-layer optimization" approach proposed in \cite{zhangSCIS}, the feedback gain of the optimal controller is
calculated by two Riccati equations independently of the Kalman filtering.
It is noted that the optimal controller at the terminal time is required to be deterministic.

The remainder of the paper is organized as follows. Section \ref{s2} presents the studied problem and some preliminaries.
The results for irregular LQ control problem are given in Section \ref{s3}.
Some concluding remarks are given in the last section.

The following notations will be used throughout this paper: $R^n$
denotes the family of $n$ dimensional vectors. $x'$ means the
transpose of $x.$ A symmetric matrix $M>0\ (\geq 0)$ means
strictly positive--definite (positive semi-definite). $Range(M)$ represents the range of the matrix $M.$
$M^{\dag}$ means the Moore-Penrose inverse of the matrix $M$.

\section{Problem Formulation}\label{s2}

In this paper, we mainly consider the system governed by the following It\^{o} stochastic differential equation:
\begin{eqnarray}
dx(t)=[A(t)x(t)+B(t)u(t)]dt+D(t)dw(t),\label{op1}
\end{eqnarray}
where $x(t)\in R^n$ is the state, $u(t)\in R^m$ is the control input, $w(t)\in R^n$ is a standard Brownian motion.
The matrices $A(t),B(t),D(t)$ are deterministic and time-varying with compatible dimensions. The initial value is given by $x(t_0)=x_0$ where $x_0$
is a random variable with mean $\hat{x}_0$ and is independent with $w(t)$.

The output of the system (\ref{op1}) is described by
\begin{eqnarray}
dy(t)=C(t)x(t)dt+G(t)dv(t),y(t_0)=0,\label{op2}
\end{eqnarray}
where $y(t)\in R^s$ is channel output, $v(t)\in R^s$ is a standard Brownian motion which is independent with $w(t)$ and $x_0$.
The matrices $C(t)$ and $G(t)$ are deterministic and time-varying with compatible dimensions.

As has been discussed in Section \ref{s1}, in order to ensure the solvability of the optimization problem, we define the following cost function:
\begin{eqnarray}
J(t_0,x_0;u)&=&E\int_{t_0}^T\Big(x'(t)Q(t)x(t)+u'(t)R(t)u(t)\Big)dt\nonumber\\
&&+[Ex'(T)]H[Ex(T)],\label{op3}
\end{eqnarray}
where $Q(t),R(t),H$ are positive semi-definite matrices with compatible dimensions. The main difference of the cost function (\ref{op3})
from that in the classic output feedback control problem lies in that the terminal term is as $[Ex'(T)]H[Ex(T)]$ rather than $E[x'(T)Hx(T)]$.

Since only partial information of the state is available to the control design, that is, $u(t)$ depends only on the
$y(s)$ for $0\leq s\leq t$, to express this non-anticipative dependence, we introduce the set
\begin{eqnarray}
\mathcal{Y}_t=\sigma\{y(s),0\leq s\leq t\}.\nonumber
\end{eqnarray}

The optimization problem studied in the paper is stated in details as follows.

\textbf{Problem (Output Feedback LQ):} For any given initial pair $(t_0,x_0),$ find a $u^*(t)$ which is $\mathcal{Y}_t$-adapted such that
\begin{eqnarray}
J^*(t_0,x_0;u^*)=\min_{u(t)\in \mathcal{Y}_t}J(t_0,x_0;u)\nonumber
\end{eqnarray}
subject to system (\ref{op1}) and (\ref{op2}) where $J(t_0,x_0;u)$ is defined by (\ref{op3}).

Since the control $u(t)$ is constrained to be $\mathcal{Y}_t$-adapted due to the fact that $x(t)$ is partially known in the form of (\ref{op2}), it is necessary to calculate
\begin{eqnarray}
\hat{x}(t)=E[x(t)|\mathcal{Y}_t].\label{op15}
\end{eqnarray}

To this end, we first state the dynamic of the optimal estimation (\ref{op15}).
By using the standard knowledge of Kalman filtering \cite{optimalfiltering}, we have
\begin{eqnarray}
d\hat{x}(t)&=&\Big(A(t)\hat{x}(t)+B(t)u(t)\Big)dt+L(t)d\nu(t),\label{op5}
\end{eqnarray}
where $d\nu(t)=dy(t)-C(t)\hat{x}(t)dt$ is the innovation process and the matrix $L(t)$ is defined by
\begin{eqnarray}
L(t)=\hat{P}(t)C'(t)\Big(G(t)G'(t)\Big)^{-1},\label{op8}
\end{eqnarray}
while $\hat{P}(t)$ satisfies the Riccati equation:
\begin{eqnarray}
\dot{\hat{P}}(t)&=&A(t)\hat{P}(t)+\hat{P}(t)A'(t)+D(t)D'(t)\nonumber\\
&&-\hat{P}(t)C'(t)\Big(G(t)G'(t)\Big)^{-1}C(t)\hat{P}(t),\nonumber
\end{eqnarray}
with the initial value $\hat{P}(t_0)=E[(x_0-\hat{x}_0)(x_0-\hat{x}_0)']$.

By applying some calculations to the innovation process, it holds that
\begin{eqnarray}
d\nu(t)=C(t)\Big(x(t)-\hat{x}(t)\Big)dt+G(t)dv(t).\nonumber
\end{eqnarray}


Denote $\tilde{x}(t)=x(t)-\hat{x}(t)$, it follows from (\ref{op1}) and (\ref{op5}) that
\begin{eqnarray}
d\tilde{x}(t)&=&A(t)\tilde{x}(t)dt+D(t)dw(t)\nonumber\\
&&-L(t)d\nu(t).\label{op9}
\end{eqnarray}

Next, we aim to reformulate the cost function (\ref{op3}) in terms of $\hat{x}(t)$ and $\tilde{x}(t)$ and present an equivalent optimization problem
to Problem (Output Feedback LQ).

Since $x(t)=\hat{x}(t)+\tilde{x}(t)$ and
$E[\hat{x}'(t)\tilde{x}(t)]=0$, the cost function (\ref{op3})
is easily rewritten as
\begin{eqnarray}
&&J(t_0,x_0;u)\nonumber\\
&=&E\int_{t_0}^T[\hat{x}'(t)Q(t)\hat{x}(t)+u'(t)R(t)u(t)]dt\nonumber\\
&&\hspace{-3mm}+E[\hat{x}(T)]'HE[\hat{x}(T)]+E\int_{t_0}^T\tilde{x}'(t)Q(t)\tilde{x}(t)dt.\label{op10}
\end{eqnarray}
Denote \begin{eqnarray}
\hat{J}(t_0,x_0;u)&=&E\int_{t_0}^T[\hat{x}'(t)Q(t)\hat{x}(t)+u'(t)R(t)u(t)]dt\nonumber\\
&&+E[\hat{x}'(T)]HE[\hat{x}(T)].\label{op7}
\end{eqnarray}
Note that the last term of (\ref{op10}) of $E\int_{t_0}^T\tilde{x}'(t)Q(t)\tilde{x}(t)dt$ is not related with the control $u$. Then, Problem (Output Feedback LQ) is reduced to the following Problem (LQG).

\textbf{Problem (LQG )}: For any given initial pair $(t_0,x_0),$ find a $\mathcal{Y}_t$-adapted $u^*(t)$ which is such that
\begin{eqnarray}
\hat{J}^*(t_0,x_0;u^*)=\min_{u(t)\in \mathcal{Y}_t}\hat{J}(t_0,x_0;u)\nonumber
\end{eqnarray}
subject to system (\ref{op5}).

\begin{theorem}\label{t1}
The optimal solution of Problem (LQG) is the optimal solution of Problem (Output Feedback LQ).
\end{theorem}
\emph{Proof.}
From (\ref{op9}), it is seen that the estimation error $\tilde{x}(t)$ is independent with the control $u(t).$
Combining with (\ref{op10}), the optimal solution of Problem (Output Feedback LQ) is equivalently to minimize (\ref{op7}). The proof is now completed.
\hfill $\blacksquare$

Based on Theorem \ref{t1}, the aim in the sequel is to solve Problem (LQG) to derive the optimal solution.

\section{Preliminaries on Irregular Optimal Control of Deterministic System}\label{s3}

Before solving Problem (LQG), we first present some preliminary results of the deterministic irregular LQ problem which have been obtained in \cite{zhangIR}.
The deterministic system is as
\begin{eqnarray}
\dot{x}(t)&=&A(t)x(t)+B(t)u(t),~x(t_0)=x_0,\label{d1}
\end{eqnarray}
and the cost function is given by
\begin{eqnarray}
J_d(t_0,x_0;u)&=&\int_{t_0}^T[x'(t)Q(t)x(t)+u'(t)R(t)u(t)]dt\nonumber\\
&&+x'(T)Hx(T),\label{d2}
\end{eqnarray}
where the matrices are the same as those in (\ref{op1}) and (\ref{op3}).

Let $P(t)$ be the solution to the following Riccati equation:
\begin{eqnarray}
0&=&\dot{P}(t)+A'(t)P(t)+P(t)A(t)+Q\nonumber\\
&&-P(t)B(t)R^{\dag}(t)B'(t)P(t), \label{dr1}
\end{eqnarray}
where the terminal condition is given by $P(T)=H$ and $R^{\dag}(t)$ represents the Moore-Penrose inverse of $R(t).$
As has been pointed in \cite{zhangSCIS}, there exist two cases with respect to $P(t)$:
\begin{eqnarray}
Range [B'(t)P(t)]\subseteq Range[R(t)],\label{d32}
\end{eqnarray}
and
\begin{eqnarray}
Range [B'(t)P(t)]\not\subseteq Range [R(t)].\label{d33}
\end{eqnarray}

In the case of (\ref{d32}), the optimal controller can be directly obtained from the equilibrium condition, which is a classic LQ optimal control problem.
In fact, the optimal controller satisfies
\begin{eqnarray}
0&=&R(t)u(t)+B'(t)p(t),\label{3}
\end{eqnarray}
where $p(t)$ satisfies
\begin{eqnarray}
\dot{p}(t)&=&-[A'(t)p(t)+Q(t)x(t)],\label{2}
\end{eqnarray}
with $p(T)=Hx(T).$ With (\ref{d32}), we can reformulate (\ref{3}) as
\begin{eqnarray}
0=R(t)u(t)+B'(t)P(t)x(t).\label{op4}
\end{eqnarray}
Thus, the linear equation (\ref{op4}) is solvable for any $x(t)\in R^n,$ which implies the solvability of the optimization problem.

While, in the case of (\ref{d33}), the linear equation (\ref{op4}) is unsolvable, that is, the control can not be solved from the standard Maximum Principle.
Nonetheless, the optimization problem may be still solvable~\cite{zhangIR}.
The optimization problem in the case of (\ref{d33}) is named as irregular LQ control~\cite{zhangSCIS}.

In this paper, we focus on the output feedback control in the irregular case (\ref{d33}). Following \cite{zhangIR}, the following denotations are firstly introduced
for convenience of future use.
Let $rank [R(t)]=m_0< m,$ thus $rank(I-R^{\dag}R)=m-m_0>0$. There is an elementary row transformation matrix $T_0(t)$ such that
\begin{eqnarray}
T_0(t)[I-R^{\dag}(t)R(t)]=\left[
                              \begin{array}{c}
                                 0  \\
                                 \Upsilon_{T_0}(t)\\
                              \end{array}
                            \right], \label{jnYC1}
\end{eqnarray}
where $\Upsilon_{T_0}(t)\in R^{[m-m_0]\times m}$ is full row rank. Furthermore denote
\begin{eqnarray}
A_0(t)&=&A(t)-B(t)R^{\dag}(t)B'(t)P(t),\nonumber\\
D_0(t)&=&-B(t)R^{\dag}(t)B'(t),\nonumber\\
 \left[ \begin{array}{cc} \ast & C_0'(t) \\
                              \end{array} \right] &=&P(t)B(t)\Big(I-R^{\dag}(t)R(t)\Big){T_0}^{-1}(t),\nonumber\\
\left[ \begin{array}{cc}\ast & B_0(t)\\
                              \end{array}
                            \right] &=&B(t)\Big(I-R^{\dag}(t)R(t)\Big){T_0}^{-1}(t),\nonumber\\
                           \left[
                                            \begin{array}{cc}
                                              * & G_0(t) \\
                                            \end{array}
                                          \right]&=& T_0(t)^{-1},\nonumber
\end{eqnarray}
and define
\begin{eqnarray}
0&=&\dot{P}_1(t)+P_1(t)A_0(t)+A_0'(t)P_1(t)\nonumber\\
&&+P_1(t)D_0(t)P_1(t), \label{d3}
\end{eqnarray}
where the terminal value $ P_1(T)$ is to be determined.

The solvability condition of the irregular LQ problem is given below.
\begin{lemma}\label{lemir}
In the case of (\ref{d33}), the control problem of minimizing the cost function $J_d(t_0,x_0;u)$ s.t. system (\ref{d1}) is solvable if and only if there exists $P_1(t)$ of (\ref{d3}) with terminal value $P_1(T)$ such that
\begin{eqnarray}
0=C_0(t)+B_0'(t)P_1(t), ~~~0\leq t\leq T,\label{d4}
\end{eqnarray}
and there exists $u_{1}(t)$ to achieve
\begin{eqnarray}
P_1(T)x(T)=0,\label{d5}
\end{eqnarray} where $x(t)$ obeys
\begin{eqnarray}
\dot{x}(t)&=&\Big(A_0(t)+D_0P_1(t)\Big)x(t)+B_0u_1(t),\label{d24}
\end{eqnarray}
with initial value $x(t_0)=x_0.$ In this case, the optimal controller $u(t)$ is given by
\begin{eqnarray}
u(t)&=&-R^{\dag}(t)B'(t)\Big(P(t)+P_1(t)\Big)x(t)\nonumber\\
&&+G_0(t)u_1(t).\nonumber
\end{eqnarray}
and the optimal cost is given by
\begin{eqnarray}
J_T^*=x_0'\Big(P(t_0)+P_1(t_0)\Big)x_0.\nonumber
\end{eqnarray}

\end{lemma}

{\em Proof}. Please see \cite{zhangIR} and \cite{zhangSCIS} for details. \hfill $\blacksquare$

\section{Solution to Irregular Problem (LQG)}\label{s4}

Having given the results for irregular LQ in last section, we are now in the position to study the irregular problem (LQG), that is, (\ref{d33}) holds. Two steps are to be presented. The first step is to
derive the FBDEs by applying the maximum principle with respect to the modified cost function (\ref{op3}).
The second step is to solve the FBDEs and obtain the optimal solution of Problem (LQG).

First, we show that the solvability of Problem (LQG) is equivalent to that of FBDEs by using the maximum principle.
\begin{lemma}\label{mp}
Problem (LQG) is solvable if and only if there exists solution $(\hat{x}(t),p(t),q(t),u(t))$ satisfying the FBDEs
\begin{eqnarray}
d\hat{x}(t)&=&\Big(A(t)\hat{x}(t)+B(t)u(t)\Big)dt+L(t)d\nu(t),\label{op16}\\
dp(t)&=&-[A'(t)p(t)+Q(t)\hat{x}(t)]dt+q(t)d\nu(t),\label{op17}\\
0&=&R(t)u(t)+B'(t)p(t),\label{op18}
\end{eqnarray}
with initial value $x(t_0)=x_0$ and terminal value $p(T)=HE[x(T)].$
\end{lemma}
\emph{Proof.}
``Sufficiency" Let $\mu(t)$ be an arbitrary control and $u(t)$ be the control satisfying (\ref{op18}). Denote $\hat{x}^\mu(t),\hat{J}(t_0,x_0;\mu)$ and $\hat{x}(t),\hat{J}(t_0,x_0;u)$
the corresponding states and cost functions with respect to $\mu(t)$ and $u(t)$ respectively. It will be shown that $\hat{J}(t_0,x_0;\mu)-\hat{J}(t_0,x_0;u)\geq0.$ In fact,
\begin{eqnarray}
&&\hat{J}(t_0,x_0;\mu)-\hat{J}(t_0,x_0;u)\nonumber\\
&=&E\int_{t_0}^T\Big([\hat{x}^{\mu}(t)-\hat{x}(t)]'Q(t)[\hat{x}^\mu(t)-\hat{x}(t)]\nonumber\\
&&+[\mu(t)-u(t)]'R(t)[\mu(t)-u(t)]\Big)dt+\Big(E[\hat{x}^\mu(T)]\nonumber\\
&&-E[x(t)]\Big)'H\Big(E[\hat{x}^\mu(T)]-E[\hat{x}(T)]\Big)\nonumber\\
&&+2E\int_{t_0}^T\Big([\hat{x}^{\mu}(t)-\hat{x}(t)]'Q(t)\hat{x}(t)\nonumber\\
&&+[\mu(t)-u(t)]'R(t)u(t)\Big)dt+2\Big(E[\hat{x}^\mu(T)]\nonumber\\
&&-E[\hat{x}(T)]\Big)'HE[x(T)]\nonumber\\
&\geq &2E\int_{t_0}^T\Big([\hat{x}^{\mu}(t)-\hat{x}(t)]'Q(t)\hat{x}(t)\nonumber\\
&&+[\mu(t)-u(t)]'R(t)u(t)\Big)dt+2\Big(E[\hat{x}^\mu(T)]\nonumber\\
&&-E[\hat{x}(T)]\Big)' HE[\hat{x}(T)].\label{op19}
\end{eqnarray}
From (\ref{op5}), it yields that
\begin{eqnarray}
d\hat{x}^\mu(t)&=&\Big(A(t)\hat{x}^\mu(t)+B(t)\mu(t)\Big)dt+L(t)d\nu(t),\nonumber\\
d\hat{x}(t)&=&\Big(A(t)\hat{x}(t)+B(t)u(t)\Big)dt+L(t)d\nu(t).\nonumber
\end{eqnarray}
This gives that
\begin{eqnarray}
&&d[\hat{x}^\mu(t)-\hat{x}(t)]\nonumber\\
&=&\Big(A(t)[\hat{x}^\mu(t)-\hat{x}(t)]+B(t)[\mu(t)-u(t)]\Big)dt.\nonumber
\end{eqnarray}
Combining with (\ref{op17}), one has
\begin{eqnarray}
&&d\Big([\hat{x}^\mu(t)-\hat{x}(t)]'p(t)\Big)\nonumber\\
&=&\Big(A(t)[\hat{x}^\mu(t)-\hat{x}(t)]+B(t)[\mu(t)-u(t)]\Big)'p(t)dt\nonumber\\
&&-[\hat{x}^\mu(t)-\hat{x}(t)]'[A'(t)p(t)+Q(t)x(t)]dt\nonumber\\
&&+[\hat{x}^\mu(t)-\hat{x}(t)]'q(t)dw(t)\nonumber\\
&=&-[\hat{x}^\mu(t)-\hat{x}(t)]'Q(t)x(t)dt\nonumber\\
&&+[\mu(t)-u(t)]'B'(t)p(t)dt\nonumber\\
&&+[\hat{x}^\mu(t)-\hat{x}(t)]'q(t)dw(t).\nonumber
\end{eqnarray}
By taking integral from $t_0$ to $T$, it is further obtained that
\begin{eqnarray}
&&E[\hat{x}^\mu(T)-\hat{x}(T)]'HE[\hat{x}(T)]\nonumber\\
&=&-E\int_{t_0}^T\Big([\hat{x}^\mu(t)-\hat{x}(t)]'Q(t)\hat{x}(t)\nonumber\\
&&+[\mu(t)-u(t)]'B'(t)p(t)\Big)dt.\nonumber
\end{eqnarray}
By substituting the above equation into (\ref{op19}), we can rewrite (\ref{op19}) as
\begin{eqnarray}
&&\hat{J}(t_0,x_0;\mu)-\hat{J}(t_0,x_0;u)\nonumber\\
&\geq &2E\int_{t_0}^T\Big([\mu(t)-u(t)]'[R(t)u(t)+B'(t)p(t)]\Big)dt\nonumber\\
&=&0,\nonumber
\end{eqnarray}
where the derivation of the last equality depends on (\ref{op18}).
This immediately gives the fact that $u(t)$ defined by (\ref{op18}) is optimal.

``Necessity" The proof is similar to that in \cite{juanPhd}, we only state the outline and omit the detailed derivations.
The main procedure is to study the variation of the cost function. To this end, let $u(t), \delta u(t)$ be the admissible controllers.
Then $v(t)=u(t)+\varepsilon \delta u(t)$ is also admissible for $\varepsilon\in[0,1]$.
Denote the corresponding states and cost functions be $\hat{x}^v(t),J(t_0,x_0;v)$ and $\hat{x}(t),J(t_0,x_0;u)$ with respect to $v(t)$ and $u(t)$ respectively.
Let $\delta \hat{x}(t)=\hat{x}^v(t)-\hat{x}(t)$ and $\delta J=J(t_0,x_0;v)-J(t_0,x_0;u),$ it is obtained following \cite{juanPhd} that
\begin{eqnarray}
\delta J=2E\int_{t_0}^T\varepsilon(\delta u(s))'\Big(R(s)u(s)+B'(s)p(s)\Big)ds+o(\varepsilon),\nonumber
\end{eqnarray}
where $p(s)$ obeys the dynamic (\ref{op17}).
This implies that the optimal controller satisfies (\ref{op18}). Thus, the solvability of FBDEs (\ref{op16})-(\ref{op18}) follows from the
solvability of Problem (LQG).
The proof is now completed. \hfill $\blacksquare$

\begin{remark}
From Lemma \ref{mp}, it is seen that the terminal value of $p(t)$ satisfies that $p(T)=HE[x(T)]$ which is deterministic. Combining this with (\ref{op18}),
we obtain that $u(T)$ is constrained to be $\mathcal{F}_0$-adapted. This is what makes Problem (LQG) solvable in the framework of the modified cost function rather than the classic
cost function.
\end{remark}

Next, we will discuss the explicitly optimal solution to Problem (LQG). The key is the solving of FBDEs (\ref{op18})-(\ref{op17}) based on Lemma \ref{mp}.
Combining with the preliminaries in Lemma \ref{lemir},
we give the result on the solvability of Problem (LQG) as shown below.

\begin{theorem}\label{t2}
In the case of (\ref{d33}), Problem (LQG) is solvable if and only if there exists $P_1(t)$ of (\ref{d3}) with terminal value $P_1(T)$ such that
(\ref{d4}) holds and there exists $u_{1}(t)$ to achieve
\begin{eqnarray}
P_1(T)E[\hat{x}(T)]=0,\label{op23}
\end{eqnarray} where $\hat{x}(t)$ obeys
\begin{eqnarray}
d\hat{x}(t)&=&\Big([A_0(t)+D_0P_1(t)]\hat{x}(t)+B_0u_1(t)\Big)dt\nonumber\\
&&+L(t)d\nu(t),\label{op24}
\end{eqnarray}
with initial value $x(t_0)=x_0.$ In this case, the optimal solution to Problem (LQG) is given by
\begin{eqnarray}
u(t)&=&\left\{
         \begin{array}{ll}
           -R^{\dag}(t)B'(t)[P(t)+P_1(t)]\hat{x}(t)\\
+G_0(t)u_1(t), & \hbox{$t<T$;} \\
           -R^{\dag}(t)B'(t)HE[x(t)]\\
+[I-R^{\dag}(t)R(t)]\xi, & \hbox{$t=T$,}
         \end{array}
       \right.\label{irreg1}
\end{eqnarray}
where $\xi$ is an arbitrary vector with compatible dimension.

\end{theorem}
\emph{Proof.} ``Necessity"
From Lemma \ref{mp}, the solvability of Problem (LQG) implies that of FBDEs (\ref{op16})-(\ref{op18}). Now we show that the necessary condition can be obtained from
the solvability of FBDEs (\ref{op16})-(\ref{op18}). We start from solving the FBDEs (\ref{op16})-(\ref{op18}).  To this end, we define
\begin{eqnarray}
\Theta(t)=p(t)-P(t)\hat{x}(t),\label{op27}
\end{eqnarray}
where $P(t)$ satisfies (\ref{dr1}) and the terminal value is given by $\Theta(T)=p(T)-HE[x(T)].$ From $p(T)=HE[x(T)],$ we have that $\Theta(T)=0.$
Assume that $d\Theta(t)=\Theta_1(t)dt+\bar{\Theta}(t)d\nu(t),$ our next aim is to determine $\Theta_1(t)$ and $\bar{\Theta}(t).$

By substituting (\ref{op27}) into (\ref{op18}), we have
\begin{eqnarray}
0&=&R(t)u(t)+B'(t)P(t)\hat{x}(t)+B'(t)\Theta(t).\nonumber
\end{eqnarray}
This gives that
\begin{eqnarray}
u(t)&=&-R^{\dag}(t)\Big(B'(t)P(t)\hat{x}(t)+B'(t)\Theta(t)\Big)\nonumber\\
&&+\Big(I-R^{\dag}(t)R(t)\Big)z(t),\label{op28}
\end{eqnarray}
and
\begin{eqnarray}
0=\Big(I-R(t)R^{\dag}(t)\Big)\Big(B'(t)P(t)\hat{x}(t)+B'(t)\Theta(t)\Big),\label{op29}
\end{eqnarray}
where $z(t)$ is an arbitrary vector with compatible dimension.

By taking It\^{o}'s formula to (\ref{op27}) and using (\ref{op28}), it follows that
\begin{eqnarray}
dp(t)&=&\dot{P}(t)\hat{x}(t)dt+P(t)\Big[A(t)\hat{x}(t)-B(t)R^{\dag}(t)\Big(B'(t)\nonumber\\
&&\times P(t)\hat{x}(t)+B'(t)\Theta(t)\Big)+B(t)\Big(I-R^{\dag}(t)\nonumber\\
&&\times R(t)\Big) z(t)\Big]dt+P(t)L(t)d\nu(t)+\Theta_1(t)dt\nonumber\\
&&+\bar{\Theta}(t)d\nu(t).\nonumber
\end{eqnarray}
From (\ref{op18}), one has
\begin{eqnarray}
dp(t)&=&-[A'(t)P(t)\hat{x}(t)+A'(t)\Theta(t)+Q(t)\hat{x}(t)]dt\nonumber\\
&&+q(t)d\nu(t).\nonumber
\end{eqnarray}
Making comparison between the above two equations yields that
\begin{eqnarray}
0&=&\dot{P}(t)\hat{x}(t)+P(t)\Big[A(t)\hat{x}(t)-B(t)R^{\dag}(t)\Big(B'(t)\nonumber\\
&&\times P(t)\hat{x}(t)+B'(t)\Theta(t)\Big)\Big]+A'(t)P(t)\hat{x}(t)\nonumber\\
&&+A'(t)\Theta(t)+Q(t)\hat{x}(t)+P(t)B(t)\Big(I-R^{\dag}(t)\nonumber\\
&&\times R(t)\Big)z(t)+\Theta_1(t),\label{op25}\\
q(t)&=&\bar{\Theta}(t)+P(t)L(t),\label{op26}
\end{eqnarray}
By using the Riccati equation (\ref{dr1}) of $P(t),$ we can reformulate (\ref{op25}) as
\begin{eqnarray}
0&=&\Big(A'(t)-P(t)B(t) R^{\dag}(t)B'(t)\Big)\Theta(t)\nonumber\\
&&+P(t)B(t)\Big(I-R^{\dag}(t)R(t)\Big)z(t)+\Theta_1(t).\nonumber
\end{eqnarray}
Recalling the denotations in Sections III, it is obtained from the above equation that
\begin{eqnarray}
-\Theta_1(t)&=&A_0'(t)\Theta(t)+C_0'(t)u_1(t).\nonumber
\end{eqnarray}
This gives the dynamic of $\Theta(t)$ as
\begin{eqnarray}
d\Theta(t)=-\Big(A_0'(t)\Theta(t)+C_0'(t)u_1(t)\Big)dt+\bar{\Theta}(t)d\nu(t).\label{op30}
\end{eqnarray}

Together with (\ref{op29}) and (\ref{op16}), the solvability of FBDEs (\ref{op16})-(\ref{op18}) is equivalently to that there exists $u_1(t)$ such that
\begin{eqnarray}
0&=&C_0(t)\hat{x}(t)+B_0'(t)\Theta(t), \Theta(T)=0,\label{op31}
\end{eqnarray}
where $\Theta(t)$ satisfies (\ref{op30}) and $\hat{x}(t)$ obeys the dynamic:
\begin{eqnarray}
d\hat{x}(t)&=&\Big(A_0(t)\hat{x}(t)+D_0(t)\Theta(t)+B_0(t)u_1(t)\Big)dt\nonumber\\
&&+L(t)d\nu(t).\label{op32}
\end{eqnarray}

Next, we will solve the FBDEs (\ref{op30})-(\ref{op32}). By applying similar procedures to \cite{zhangSCIS}, we have that there exists a homogeneous
relationship between $\Theta(t)$ and $\hat{x}(t)$ as
\begin{eqnarray}
\Theta(t)=P_1(t)\hat{x}(t),\label{op33}
\end{eqnarray} where $P_1(t)$ satisfies (\ref{d3}) and the terminal value is given by $\Theta(T)=P_1(T)E[\hat{x}(T)]$. Since $\Theta(T)=0,$ then (\ref{op23}) follows directly.
In addition, based on (\ref{op33}) and (\ref{op31}), the condition (\ref{d4}) holds.

``Sufficiency" Based on Lemma \ref{mp}, it is sufficient to prove that the FBDEs (\ref{op16})-(\ref{op18}) are solvable. To do this, we firstly verify that
\begin{eqnarray}
p(t)&=&P(t)\hat{x}(t)+P_1(t)\hat{x}(t),\label{op20}\\
q(t)&=&[P(t)+P_1(t)]L(t),\label{op21}
\end{eqnarray} where the terminal value is given by $p(T)=HE[\hat{x}(T)]$ and $\hat{x}(t)$ satisfies
\begin{eqnarray}
d\hat{x}(t)&=&\Big[\Big(A(t)-B(t)R^{\dag}(t)B'(t)[P(t)+P_1(t)]\Big)\hat{x}(t)\nonumber\\
&&+B(t)G_0(t)u_1(t)\Big)dt+L(t)d\nu(t),\nonumber
\end{eqnarray}
is the solution to (\ref{op17}).
By using It\^{o}'s formula to $P(t)\hat{x}(t)+P_1(t)\hat{x}(t),$ it yields that
\begin{eqnarray}
&&d[P(t)\hat{x}(t)+P_1(t)\hat{x}(t)]\nonumber\\
&=&[\dot{P}(t)+\dot{P}_1(t)]\hat{x}(t)dt+[P(t)+P_1(t)]\Big(A(t)\nonumber\\
&&-B(t)R^{\dag}(t)B'(t)[P(t)+P_1(t)]\Big)\hat{x}(t)dt\nonumber\\
&&+[P(t)+P_1(t)]B(t)G_0(t)u_1(t)dt\nonumber\\
&&+[P(t)+P_1(t)]L(t)d\nu(t)\label{1}\\
&=&\Big(-A'(t)[P(t)+P_1(t)]\hat{x}(t)-Q(t)\hat{x}(t)\Big)dt\nonumber\\
&&+[P(t)+P_1(t)]L(t)d\nu(t).\nonumber
\end{eqnarray}
Since $\Upsilon_{T_0}(t)$ is of full row rank, there exits a vector $z(t)\in R^m$ such that $\Upsilon_{T_0}(t)z(t)=u_1(t)$.
This further gives that
\begin{eqnarray}
&&[P(t)+P_1(t)]B(t)G_0(t)u_1(t)\nonumber\\
&=&[P(t)+P_1(t)]B(t)G_0(t)\Upsilon_{T_0}(t)z(t)\nonumber\\
&=&[P(t)+P_1(t)]B(t)T_0^{-1}(t)T_0(t)\Big(I-R^{\dag}(t)R(t)\Big)z(t)\nonumber\\
&=&[P(t)+P_1(t)]B(t)\Big(I-R^{\dag}(t)R(t)\Big)T_0^{-1}(t)\nonumber\\
&&\times T_0(t)\Big(I-R^{\dag}(t)R(t)\Big)z(t)\nonumber\\
&=&[C_0'(t)+P_1(t)B_0(t)]T_0(t)\Big(I-R^{\dag}(t)R(t)\Big)z(t)\nonumber\\
&=&0,\nonumber
\end{eqnarray}
where we have used (\ref{d4}) in the derivation of the last equality.
Thus, (\ref{1}) becomes
\begin{eqnarray}
&&d[P(t)\hat{x}(t)+P_1(t)\hat{x}(t)]\nonumber\\
&=&\Big(-A'(t)[P(t)+P_1(t)]\hat{x}(t)-Q(t)\hat{x}(t)\Big)dt\nonumber\\
&&+[P(t)+P_1(t)]L(t)d\nu(t).\nonumber
\end{eqnarray}
Together with (\ref{op17}), we have that $(p(t),q(t))$ defined by (\ref{op20}) and (\ref{op21}) is the solution
to (\ref{op17}).

Next we derive the optimal controller (\ref{irreg1}). By substituting (\ref{op20}) into (\ref{op18}), it is obtained for $t<T$ that
\begin{eqnarray}
0&=&R(t)u(t)+B'(t)[P(t)+P_1(t)]\hat{x}(t).\nonumber
\end{eqnarray}
This gives that
\begin{eqnarray}
u(t)&=&-R^{\dag}(t)B'(t)[P(t)+P_1(t)]\hat{x}(t)\nonumber\\
&&+[I-R^{\dag}(t)R(t)]z(t),~t<T,\label{op22}
\end{eqnarray}
where $z(t)$ is an arbitrary vector with compatible dimension. It is noted that
\begin{eqnarray}
&&[I-R^{\dag}(t)R(t)]z(t)\nonumber\\
&=&T_0^{-1}(t)T_0(t)[I-R^{\dag}(t)R(t)]z(t)\nonumber\\
&=&G_0(t)\Upsilon_{T_0}(t)z(t).\nonumber
\end{eqnarray}
Together with the fact that $\Upsilon_{T_0}(t)$ is of full row rank, there exists a $u_1(t)$ such that
$[I-R^{\dag}(t)R(t)]z(t)=G_0(t)u_1(t).$ Thus the case of $t<T$ in (\ref{irreg1}) follows from (\ref{op22}).
In addition, in view of $p(T)=HE[x(T)]$ and (\ref{op18}), the case of $t=T$ in (\ref{irreg1}) is obtained.
The proof is now completed.
\hfill $\blacksquare$

\begin{remark}
As is shown in Theorem \ref{t2}, the explicitly optimal controllers (\ref{irreg1}) at terminal time $T$ should be
deterministic. This is the difference of  this paper from the classic output feedback control.
\end{remark}

Finally we study the open-loop and closed-loop solvability of Problem (LQG). The detailed results are given below.
\begin{theorem}\label{Th-1}
Problem (LQG) is open-loop solvable if and only if there exists $P_1(t)$ of (\ref{d3}) such that (\ref{d4}) holds and
\begin{eqnarray}
Range \big[P_1(t_0)\big]\subseteq
Range \Big(G_1[t_0,T]\Big),\label{df4}
\end{eqnarray}
where the Gramian matrix $G_1[t_0,T]$ is defined by
\begin{eqnarray}
G_1[t_0,T]=\int_{t_0}^TP_2(t_0,s)C_0'(s)C_0(s)P_2'(t_0,s)ds,
\end{eqnarray}
while $P_2(t,s)$ satisfies
$\dot{P}_2(t,s)=-A_0'(t)P_2(t,s),~P_2(t,t)=I.$
In this case, the open-loop solution can be given by (\ref{irreg1}) while $u_1(t)$ is given by
\begin{eqnarray}
u_1(t)&=&C_0(t)P_2'(t_0,t)G_1^{\dag}[t_0,T]P_1(t_0)x_0.\label{df5}
\end{eqnarray}
\end{theorem}
\emph{Proof.}  From Theorem \ref{t2}, Problem (LQG) is solvable if and only if there exists $P_1(t)$ of (\ref{d3}) such that
(\ref{d4}) holds and there exists $u_{1}(t)$ to achieve (\ref{op23}). Thus, it is equivalent to show that the existence of $u_1(t)$ to achieve (\ref{op23}) if and only if
(\ref{df4}) holds. From the dynamic (\ref{op24}) and the open-loop property of $u_1(t)$, it is obtained that
\begin{eqnarray}
\frac{dE[\hat{x}(t)]}{dt}&=&[A_0(t)+D_0P_1(t)]E[\hat{x}(t)]+B_0u_1(t).\nonumber
\end{eqnarray}
By applying the similar discussions to Theorem 3 in \cite{zhangSCIS}, the result follows. The proof is now completed.
\hfill $\blacksquare$

We further consider the closed-loop solution. Following \cite{zhangSCIS}, we make the following denotations
$P_1(t)=\left[
         \begin{array}{cc}
           P^1_{11}(t) & P^1_{12}(t) \\
           {P^1_{12}}'(t) & P^1_{22}(t) \\
         \end{array}
       \right].$
Assume that $P_1(t)$ is singular, then there exists $\mathcal{T}_1(t)$ such that
$\mathcal{T}_1'(t)P_1(t)\mathcal{T}_1(t)=\left[
         \begin{array}{cc}
           \hat{P}(t) & 0 \\
           0 & 0 \\
         \end{array}
       \right],$
where $\hat{P}(t)$ is invertible. Let $
\dot{\mathcal{T}}_1'(t)\mathcal{T}_1(t)=\left[
      \begin{array}{c}
        \tilde{T}_{1}(t) \\
        \tilde{T}_{2}(t) \\
      \end{array}
    \right], \mathcal{T}_1'(t)B_0\mathcal{T}_1(t)=\left[
      \begin{array}{c}
        B_{1}(t) \\
        B_{2}(t) \\
      \end{array}
    \right]$
and
$\mathcal{T}_1'(t)\Big(A_0(t)+D_0P_1(t)\Big)\mathcal{T}_1(t)=\left[
                                 \begin{array}{c}
                                   \hat{A}_{1}(t) \\
                                   \hat{A}_{2}(t) \\
                                 \end{array}
                               \right].$

\begin{theorem}\label{thd3}
Assume that (\ref{d4}) holds and there exists $K(t)$ such that
\begin{eqnarray}
\tilde{T}_{1}(t)+\hat{A}_{1}(t)+B_{1}(t)\mathcal{T}_1'(t)K(t)\mathcal{T}_1(t)=\left[
     \begin{array}{cc}
       \frac{I}{t-T} & 0 \\
     \end{array}
   \right],\label{df7}
\end{eqnarray}
then Problem (LQG) is closed-loop solvable.
In this case, the closed-loop solution is given by (\ref{irreg1}) where
$u_1(t)=K(t)x(t),$
where $K(t)$ satisfies that (\ref{df7}).
\end{theorem}
\emph{Proof.} Denote
\begin{eqnarray}
y(t)=\mathcal{T}_1'(t)E[\hat{x}(t)]=\left[
                           \begin{array}{c}
                             y_1(t) \\
                             y_2(t) \\
                           \end{array}
                         \right]. \label{Cls-1}
\end{eqnarray}
By using (\ref{op24}) and the feedback control $u_1(t)=K(t)\hat{x}(t)$, we have the dynamic
\begin{eqnarray}
\dot{y}(t)
&=&\left[
     \begin{array}{c}
       \tilde{T}_{1}(t)+\hat{A}_{1}(t)+B_{1}(t)\mathcal{T}_1'(t)K(t)\mathcal{T}_1(t)\\
       \tilde{T}_{2}(t)+\hat{A}_{2}(t)+B_{2}(t)\mathcal{T}_1'(t)K(t)\mathcal{T}_1(t)\\
     \end{array}
   \right]y(t).\nonumber
\end{eqnarray}
Then by applying similar procedures to Theorem 4 in \cite{zhangSCIS},
the result follows.
The proof is now completed. \hfill $\blacksquare$

\begin{remark}
Based on Theorem \ref{t1}, the optimal solution of Problem (Output Feedback LQ) is the same as that of Problem (LQG). Thus, the optimal solution of Problem (Output Feedback LQ)
can be given as those in Theorem \ref{t2}, \ref{Th-1} and \ref{thd3}.
\end{remark}

\section{Conclusions}

In this paper, we study the irregular output feedback LQ control problem.
By introducing a modified cost function, the separation principle was shown to be hold and the explicitly optimal controller was given in the feedback form of the Kalman filtering.
In particular, the feedback gain has been calculated by two Riccati equations independently of the Kalman filtering.
It is emphasized that the optimal controller at the terminal time is required to be deterministic.



\begin{thebibliography}{0}

\bibitem{anderson} Anderson B D O, Moore J B. Optimal control: linear quadratic methods.
Englewood Cliffs, NJ: Prentice Hall, 1990.

\bibitem{optimalfiltering} Anderson B D O, Moore J B. Optimal filtering, Englewood Cliffs, N.J.: Prentice Hall, 1979.

\bibitem{Bell} Bell D J. Singular problems in optimal control-a survey. International Journal of Control, 1975, 21(2): 319-331.

\bibitem{bellman} Bellman R. The theory of dynamic programming. \emph{Bulletin of the American Mathematical Society}, 60(6): 503-516, 1954.

\bibitem{Bellman} Bellman R, Glicksberg I, Gross O. Some aspects of the mathematical theory of control processes. Rand Corporation, R-313, 1958.

\bibitem{Bonnans} Bonnans J F, Silva F J. First and second order necessary conditions for
stochastic optimal control problems. Applied Mathematics \& Optimization, 2012, 65: 403-439.


\bibitem{chenhanfu} Chen H-F. Unified controls applicable to general case under quadratic index. Acta Mathematicae Applicatae Sinica, 1982, 5(1): 45-52.

\bibitem{chenlizhou} Chen S, Li X, Zhou X. Stochastic linear quadratic regulators with indefinite control weight costs. SIAM Journal on Control and Optimization, 1998, 36(5): 1685-1702.

\bibitem{clementsanderson} Clements D, Anderson B. Singular optimal control: The linear-quadratic problem. Springer-Verlag, New York, 1978.

\bibitem{Gurman} Gurman V. The method of multiple maxima and optimization problems for space maneuvers. Proc. Second Readings of K. E. Tsiolkovskii, Moscow, 1968, 39-51.

\bibitem{Gabasov} Gabasov R, Kirillova F M. High order necessary conditions for optimality. SIAM J. Control, 1972, 10: 127-168.

\bibitem{Ho} Ho Y. Linear stochastic singular control problems. Journal of Optimization Theory and Application, 1972, 9(1): 24-31.

\bibitem{Hoehener} Hoehener D. Variational approach to second-order optimality conditions for control
problems with pure state constraints. SIAM Journal of Control and Optimization, 2012, 50: 1139-1173.

\bibitem{Hsia} Hsia T. On the existence and synthesis of optimal singular control with
quadratic performance index. IEEE Transactions on Automatic Control, 1967, 12(6): 778-779.

\bibitem{Kalman} Kalman R E. Contributions to the theory of optimal control, Bol. Soc., Mat. Mexicana, 1960, 5: 102-119.

\bibitem{Kliger} Kliger I. Discussion on the stability of the singular trajectory with respect to ``Bang-Bang" control. IEEE Transactions on Automatic Control, 1964, 9(4): 583-585.

\bibitem{Krener}  Krener A J. The high order maximal principle and its application to singular extremals. SIAM Journal of Control and Optimization, 1977, 15: 256-293.

\bibitem{Letov} Letov A M. The analytical design of control systems. Automat. Remote Control, 1961, 22: 363-372.

\bibitem{lewis} Lewis F L, Vrabie D L, Syrmos V L. Optimal control.
John Wiley \& Sons, Inc., 2012.


\bibitem{Moore} Moore J. The singular solutions to a singular quadratic minimization problem. International Journal of Control, 1974, 20(3): 383-393.


\bibitem{pinv} Penrose R. A generalized inverse of matrices, Mathematical Proceedings of the Cambridge Philosophical Society, 1955, 52: 17-19.

\bibitem{mp} Pontryagin L S, Boltyanskii V G, Gamkrelidze R V, Mishchenko E F. The mathematical theory of optimal process.
English translation. Interscience, 1962.

\bibitem{SCIS2}	Qi Q, Zhang H. Time-inconsistent stochastic linear quadratic control for discrete-time systems. SCIENCE CHINA Information Sciences, 2017, 60:120204.

\bibitem{Speyer} Speyer J, Jacobson D. Necessary and sufficient conditions for optimality for singular control problems: a limit approach. Journal of Mathematical Analysis and Applications, 1971, 34(2): 239-266.

\bibitem{sunliyong} Sun J, Li X, Yong J. Open-loop and closed-loop solvabilities for stochastic linear quadratic optimal control problems. SIAM Journal on Control and Optimization, 2016, 54(5): 2274-2308.

\bibitem{Williems} Willems J C, Kitapci A, Silverman L M. Singular optimal
control: a geometric approach. IAM Journal of Control and Optimization, 1986, 24(2): 323-337.

\bibitem{juanPhd} Xu J. Stochastic control for systems with time delay and its applications, Dissertation for Doctoral Degree, Shandong University, 2013.

\bibitem{SCIS1} Xu J, Shi J, Zhang H. A leader-follower stochastic linear quadratic differential game with time delay. SCIENCE CHINA Information Sciences, 2018, 61:112202.

\bibitem{yong} Yong J, Zhou X. Stochatic controls: Hamiltonian systems and HJB equations. Springer Verlag, 1999.

\bibitem{hszhang} Zhang H, Lin L, Xu J, Fu M. Linear quadratic regulation and stabilization of
discrete-time Systems with delay and multiplicative noise. IEEE Transactions on Automatic Control, 2015, 60(10): 2599-2613.

\bibitem{zhangIR} Zhang H, Xu J, On Irregular Linear Quadratic Control: Deterministic Case, \emph{arXiv:1711.09213}.

\bibitem{hszhang1} Zhang H, Xu J, Control for It\^{o} stochastic systems with input delay, \emph{IEEE Trans.\ Autom.\ Control}, 62(1): 350-365, 2017.

\bibitem{zhangSCIS} Zhang H, Xu J, Optimal Control with Irregular Performance, SCIENCE CHINA Information Sciences, DOI: 10.1007/s11432-018-9685-8.

\bibitem{XuZhang} Zhang H, Zhang X. Pointwise second-order necessary conditions for stochastic optimal controls, Part I: The case of convex control constraint.
SIAM Journal of Control and Optimization, 2015, 53(4): 2267-2296.

\bibitem{ZhangFangfang} Zhang F, Zhang H, Tan C, Wang W, et al. A new approach to distributed control for multi-agent systems based on approximate upper and lower bounds.
 International Journal of Control, Automation and Systems, 2017, 15(6): 2507-2515.

\end{thebibliography}
\end{document}